\documentclass[12pt]{amsart}
\newcommand{\g}{\mathfrak}
\newcommand{\ccal}{\mathcal}

\newcommand{\gtg}{\g g}

\newcommand{\hgtg}{\hat{\gtg}}
\newcommand{\tgtg}{\tilde{\gtg}}

\newcommand{\gtsl}{{\g sl}}

\newcommand{\gth}{{\g h}}
\newcommand{\nc}{{\bf C}}
\newcommand{\nq}{{\bf Q}}
\newcommand{\tens}{\dot{\otimes}}
\newcommand{\hg}{\hat{\Gamma}}

\newcommand{\nz} {{\bf Z}}

\newcommand{\cp}{{\bf CP^{1}}}

\newcommand{\co}{{\ccal O}}
\newcommand{\tco}{{\tilde{\co_{k}}}}

\newcommand{\cf}{{\ccal F}}

\newcommand{\ca}{{\ccal A}}

\newtheorem{theorem}{Theorem}[subsection]
\newtheorem{proposition}[theorem]{Proposition}

\newtheorem{lemma}[theorem]{Lemma}
\newtheorem{corollary}[theorem]{Corollary}

\newtheorem{conjecture}[theorem]{Conjecture}

\title[Kazhdan-Lusztig tensoring and Harish-Chandra Categories]
{Kazhdan-Lusztig Tensoring and
 Harish-Chandra Categories}
\author{Igor B. Frenkel}
\address{Department of Mathematics, Yale University}
\author{Feodor Malikov}
\address{Department of Mathematics, University of Southern
California}

\normalsize 
\divide\textheight by \baselineskip
\multiply\textheight by \baselineskip
\begin{document}

\begin{abstract}

We use the Kazhdan-Lusztig tensoring to define affine translation
functors, describe annihilating ideals of highest weight modules
over an affine Lie algebra in terms of the corresponding VOA, and
to sketch a functorial approach to ``affine Harish-Chandra
bimodules''.

\end{abstract}

\maketitle

 \section{introduction}

Our original motivation was to answer the question: ``What is a Harish-Chandra
bimodule over an affine Lie algebra?'' Although we have not yet been able
to give a complete answer, we can state a conjecture and we can produce
objects which are remarkably (and non-trivially) reminiscent of the
principal series representations of a complex group. Along the way we get
a couple of results (on annihilating ideals of highest
weight modules, and on equivalence of categories) which are apparently
interesting by themselves.

\subsection{Representations of complex groups. }
\label{reviewofbernstgelf}
To make things clearer, we first review a categorical approach to
Harish-Chandra
bimodules over a simple Lie algebra following the beautiful paper by
Bernstein and S.Gelfand \cite{bernst_gelf}.  Let $\ca$ be a category.
 Then
one can consider the category $Funct(\ca)$ of functors on $\ca$,
objects
being functors, morphisms being natural transformations of functors.
In general,  there is
no reason to think that $Funct(\ca)$ is abelian even if
$\ca$ is so. Here is, however, an important example when $Funct(\ca)$ contains
an abelian complete subcategory.

Let $Mod(\gtg)$ be the category of  modules over a simple complex Lie algebra
$\gtg$ and $Mod(\gtg-\gtg)$ the category of $\gtg$-bimodules.
``Module'' will always mean a space carrying a left action of $\gtg$;
``bimodule'' will always mean a space carrying a left and a right action
commuting with each other.
Any $H\in  Mod
(\gtg-\gtg)$ gives rise to the functor
\[\Phi_{H}: Mod(\gtg)\rightarrow Mod(\gtg);\; \Phi_{H}(M)=H\otimes_{\gtg}M.\]

It is well-known that
\[Hom_{Mod(\gtg-\gtg)}(H_{1},H_{2})=Hom_{Funct(Mod(\gtg))}(\Phi_{H_{1}},
\Phi_{H_{2}}).\]
\begin{sloppypar}
Therefore $Mod(\gtg-\gtg)$ is a complete abelian subcategory of
$Funct(Mod(\gtg))$.
\end{sloppypar}

Any $\gtg$-bimodule is a $\gtg$-module with respect to the diagonal action
(that is, the left action minus the right action). A Harish-Chandra module is
a bimodule such that under the diagonal action it decomposes in a direct
sum of finite dimensional $\gtg$-modules occurring with finite multiplicities.
Consider the category of Harish-Chandra bimodules $HCh$, and the $\co$ category
of  $\gtg$-modules. The condition imposed
on the diagonal action ensures that if $H$ is a Harish-Chandra
bimodule, then $\Phi_{H}$ preserves $\co$.
Therefore the construction we just discussed gives an
embedding $HCh\hookrightarrow Funct(\co)$ as a complete subcategory.

Further, indecomposable
projective Harish-Chandra bimodules are exactly those corresponding
to direct indecomposable summands of the functor of tensoring by
a finite dimensional $\gtg$-module $V$:
\begin{equation}
\label{projfunctdef}
V\otimes ?: \co\rightarrow\co,\; M\mapsto V\otimes M.
\end{equation}
Such functors are called {\em projective}.

Having classified projective functors, it is relatively easy to establish an
equivalence (see (\ref{bgisocat1}) below) of (sub)categories of $HCh$ and
$\co$.

To be more precise, observe that
  $HCh$ and $\co$ admit direct product decompositions with
respect to the
action of the center of the universal enveloping $U(\gtg)$.  Namely,
\[ HCh=\oplus_{\theta_{l},\theta_{r}}HCh(\theta_{l},\theta_{r}),
\co=\oplus_{\theta}\co_{\theta},\]
where $\theta_{l},\theta_{r},\theta$ are central characters,
$HCh(\theta_{l},\theta_{r})\subset HCh$ is a complete subcategory of
 Harish-Chandra bimodules admitting left central character $\theta_{l}$ and
right central character $\theta_{r}$; $\co_{\theta}\subset\co$ is a complete
subcategory defined in a similar way.  It is  easy to see that
$HCh(\theta_{l},\theta_{r})$ is empty unless $\lambda_{l}-\lambda_{r}$ is
integral, where $\lambda_{l}$ (resp. $\lambda_{r}$) is a dominant weight
related to
$\theta_{l}$ (resp. $\theta_{r}$); this condition will be tacitly assumed
from now no.

Another convention to be adopted for simplicity
is that all
central characters in question are assumed to be regular, i.e. all
corresponding weights
are off the walls of the Weyl chambers.

One of the main results of \cite{bernst_gelf} is that the functor

\begin{equation}
\label{bgisocat1}
HCh(\theta_{l},\theta_{r})\rightarrow \co_{\theta_{l}},\;
H\mapsto \Phi_{H}(M_{\lambda_{r}}),
\end{equation}
is an equivalence of categories. Here $M_{\lambda_{r}}$ is the Verma module
with the highest dominant weight $\lambda_{r}$.

\begin{sloppypar}
{}From this one gets the principal series representations
$H_{w}\in HCh(\theta_{l},\theta_{r}) ,w\in W$, as
preimages
of the Verma modules  $M_{w\cdot\lambda_{l}}$ under (\ref{bgisocat1}):
\[H_{w}=Hom_{\nc}(M_{\lambda_{r}}, M_{w\cdot\lambda_{l}})^{fin},\]
where $Hom_{\nc}(M_{\lambda_{r}}, M_{w\cdot\lambda_{l}})$ is understood
as a $\gtg$-bimodule with respect to the obvious bimodule structure
and $Hom_{\nc}(M_{\lambda_{r}}, M_{w\cdot\lambda_{l}})^{fin}\subset
Hom_{\nc}(M_{\lambda_{r}}, M_{w\cdot\lambda_{l}})$ is the maximal submodule
locally finite with respect to the diagonal action. Thus we get, in particular,
that
simple Harish-Chandra bimodules are labelled by the elements of the Weyl
group $W$.
\end{sloppypar}

Another important corollary of (\ref{bgisocat1}) is the following
description of the 2-sided ideal lattice of $U(\gtg)_{\theta}:=
U(\gtg)/U(\gtg)Ker(\theta)$. Denote by $\Omega(U(\gtg)_{\theta})$ the
2-sided ideal lattice of $U(\gtg)_{\theta}$ and by $\Omega(M_{\lambda})$
the submodule lattice of $M_{\lambda}$, where $\lambda$ is the dominant
weight related to $\theta$. Then
the map

\begin{equation}
\label{bg_descrofid}
\Omega(U(\gtg)_{\theta})\rightarrow \Omega(M_{\lambda}),
I\mapsto  IM_{\lambda}
\end{equation}
is a lattice equivalence. Indeed,  $U(\gtg)_{\theta}$
is an algebra containing $\gtg$, and hence a $\gtg$-bimodule; its 2-sided
ideals as algebra
 are its
submodules as bimodule. Under the equivalence (\ref{bgisocat1})
$U(\gtg)_{\theta}$ goes to $M(\lambda)$, because
$U(\gtg)_{\theta}\otimes_{\gtg}M(\lambda)=M(\lambda)$. Thus submodule lattices
of $U(\gtg)_{\theta}$ and $M(\lambda)$ are equivalent. A little extra work
is needed to find the explicit form (\ref{bg_descrofid}) of this
equivalence.

 The last result of \cite{bernst_gelf} which we want to review here is
another equivalence of categories based on the notion of a translation
functor. Let the central characters $\theta_{1},\theta_{2}$ be such that
the difference of the corresponding dominant highest weights $\lambda_{1}-
\lambda_{2}$ is integral. Denote by $\lambda$ the dominant weight lying
in the $W$-orbit of  $\lambda_{1}-
\lambda_{2}$, and by $V_{\lambda}$ the simple $\gtg$-module with highest
weight $\lambda$.
 For any $\theta$ denote by $p_{\theta}:\co\rightarrow\co_{\theta}$
the natural projection. Then the functor

\begin{equation}
\label{bf_equiv2}
T_{\theta_{2}}^{\theta_{1}}:\co_{\theta_{2}}\rightarrow \co_{\theta_{1}},\;
T_{\theta_{2}}^{\theta_{1}}(M)=p_{\theta_{1}}(V_{\lambda}\otimes M)
\end{equation}
is an equivalence of categories. The functor $T_{\theta_{2}}^{\theta_{1}}$
is called translation functor.

We finish our review of the semi-simple case by remarking that many results
of \cite{bernst_gelf} are based on, refine and generalize the earlier work,
see e.g. \cite{enr,duflo,vogan,zuck}.

\subsection{An affine analogue.}
There are many reasons why it is difficult to give
an intelligent definition of a Harish-Chandra bimodule over an affine
Lie algebra $\hgtg$. (For one thing, it  follows from our
results that one should rather define a Harish-Chandra bimodule over
the corresponding vertex operator algebra.) We find it easiest to
adopt a functorial point of view.

Thus we are looking for an interesting subcategory in $Funct(\tco)$,
$\tco$ being the Bernstein-Gelfand-Gelfand category of $\hgtg$-modules
at level $k$
satisfying the additional condition that the modules are semi-simple
over $\gtg\subset\hgtg$. As an analogue  of the functor $V\otimes ?$
we choose
\[V_{\lambda}^{k}\tens ?: \tco\rightarrow\tco,\; A\mapsto V_{\lambda}^{k}\tens
A,\]
where $\tens:\tco\times\tco\rightarrow\tco$ is the Kazhdan-Lusztig tensoring
\cite{kazh_luszt_0,kazh_luszt,kazh_luszt_1},
and $V_{\lambda}^{k}$ is the Weyl module (generalized Verma module in another
termiology) induced in a standard way from the finite dimensional $\gtg$-module
$V_{\lambda}$. (There seems to be no other reasonable choice.)

The Kazhdan-Lusztig tensoring is a subtle thing and many obvious properties of
$V\otimes ?$ are hard to carry over to the case of $V_{\lambda}^{k}\tens ?$.
For example, the functor $V_{\lambda}^{k}\tens ?$ does not seem to be
exact in general.
There is, however, a case when the analogy is precise -- the affine version
of a translation functor. By \cite{deodgabbkac,rcw}, there is a direct sum
decomposition
\[\tco=\oplus_{(\lambda,k)\in P^{+}_{k}}\tco^{\lambda},\]
and thus a projection
\[p_{\lambda}:\tco\rightarrow\tco^{\lambda},\]
where $P^{+}_{k}$ is the set of dominant weights at level
$k+h^{\vee}\in\nq_{>}$. (This is
an analogue of the central character decomposition for $\gtg$.) We can
therefore
define an affine translation functor
\[T_{\mu}^{\lambda}:\tco^{\mu}\rightarrow \tco^{\lambda},\]
by adjusting definition (\ref{bf_equiv2}) to the affine case (most notably by
replacing
$\otimes$ with $\tens$ and the finite dimensional $\gtg$-module with an
appropriate
Weyl module, for details see \ref {defoftrfunctor}). This construction
was first proposed in \cite{fink} in the case of  negative level
($k+h^{\vee}< 0$)
representations.

The basic properties of affine translation functors are collected
in Proposition \ref{conseqreltochar}. They are summarized by
saying that a Weyl module with a dominant highest weight is rigid
and the functor of Kazhdan-Lusztig tensoring with such a module
is exact.  These prperties easily imply that
$T_{\mu}^{\lambda}:\tco^{\mu}\rightarrow\tco^{\lambda}$ is
an equivalence of categories (c.f. (\ref{bf_equiv2}). This
theorem refines
results of \cite{deodgabbkac}, where a different version
of translation functors was
defined (in the framework of a general symmetrizable Kac-Moody
algebra) by using the standard tensoring with an integrable module.

The study of Kazhdan-Lusztig tensoring is not easy but rewarding.
A simple translation of Proposition \ref{conseqreltochar} in
the language of vertex operator algebras
(see \ref{Vertexoperatorsandvertexoperatoralgebras},
\ref{and_vertex_operator_algebras})
gives the following
affine analogue of the equivalence (\ref{bg_descrofid}). Recall
that by \cite{frzhu}, there is a vertex operator algebra (VOA)
$(V_{0}^{k},Y(.,t))$ attached to  $\hgtg$. The Fourier components
of the fields $Y(v,t),v\in V_{0}^{k}$ span a Lie algebra,
$U(\hgtg)_{loc}$. We prove (Theorem \ref{theoronannideals}) that
the ideal lattice of $U(\hgtg)_{loc}$ as VOA is equivalent to
the submodule lattice of the Weyl module $V_{\lambda}^{k}$
with a dominant highest weight $(\lambda,k),\;k+h^{\vee}\in\nq_{>}$.
Observe that the crucial difference between this statement
and (\ref{bg_descrofid}) is that the asociative algebra
$U(\gtg)_{\theta}$ is replaced by a huge Lie algebra
 $U(\hgtg)_{loc}$. Theorem \ref{theoronannideals} generalizes
and refines the well-known result that Fourier components of
the field $e_{\theta}(t)^{k+1}$ annihilate all integrable
modules at a positive level $k$; here $e_{\theta}\in\gtg$ is
a highest root vector.

Having found two affine analogues of two corollaries of the
equivalence (\ref{bgisocat1}), we return to the problem
of affinizing the notion of a Harish-Chandra module. We conjecture
(for details see sect.\ref{whcafm}) that the functor

\[\tco^{0}\rightarrow Funct(\tco^{\lambda}),\; A\mapsto p_{\lambda}
\circ(A\tens ?)\]

realizes $\tco^{0}$ as a complete subcategory of
$Funct(\tco^{\lambda})$. Realized in this way $\tco^{0}$
becomes a precise analogue of $HCh(\lambda,\lambda)$.
(We are forced to change notation from  $HCh(\theta,\theta)$
to  $HCh(\lambda,\lambda)$ as introducing the notion
of a central character is troublesome in the affine case.) The theorem
on affine translation functors then shows that (affine)
$HCh(\lambda,\lambda)$ is equivalent to $\tco^{\lambda}$ as it
should in light of (\ref{bgisocat1}). As
a supporting evidence we show that the natural map

\[Hom_{\hgtg}(A, V_{w\cdot 0}^{k})\rightarrow
Hom_{\hgtg}(p_{\lambda}(A\tens V_{\lambda}^{k}), p_{\lambda}
(V_{w\cdot 0}^{k}\tens
V_{\lambda}^{k}))\]

is an isomorphism. Therefore
 there is an injection
\[Hom_{\hgtg}(A, V_{w\cdot 0}^{k})\hookrightarrow
Hom_{Funct(\tco^{\lambda})}(p_{\lambda}(A\tens ?),
p_{\lambda}(V_{w\cdot 0}^{k}\tens ?)).\]
(The conjecture would imply that this map is an isomorphism.)
Thus the functor
 \[p_{\lambda}\circ(V_{w\cdot 0}^{k}\tens ?)
:\tco^{\lambda}\rightarrow \tco^{\lambda},\]
is indeed very reminiscent of the principal series
 representation $H_{w}$, insofar as the Weyl module $V_{w\cdot
0}^{k}$ is analogous to the Verma module $M_{w\cdot 0}$.

We finally observe that all these have analogues for the category
$\co_{k}\supset\tco$ obtained by dropping the condition of
$\gtg$-semi-simplicity. In this way we get objects better modelling
principal series representations in the affine case.

\begin{sloppypar}
{\bf Acknowledgements.} We thank P.Etingof, B.Feigin and
G.Zuckerman for
illuminating discussions.
Part of the work was done when we were visiting E.Shrodinger Institute
in Vienna in June 1996 . We are grateful to the organizers of the conference
on Representation  Theory and Applications to Mathematical Physics for
invitation and to the Institute for support. The second author enjoyed
the support of IHES and discussions with K.Gawedzki during July-August of
1996.
\end{sloppypar}

\section{preliminaries}
\label{preliminaries}
\subsection{ }
\label{prelim_reps}
  The following is a list of essentials which
will be used but will not be explained.

$\gtg$ is a simple finite dimensional Lie algebra
with a fixed triangular decomposition;
in particular with a fixed Cartan subalgebra
$\gth\subset\gtg$;

the action ($\lambda\mapsto
w\lambda$) and the shifted (by $\rho$)
 action ($\lambda\mapsto w\cdot\lambda$) action
of the Weyl group $W$ on $\gth^{\ast}$ preserving the
weight lattice $P\in\gth^{\ast}$; denote by $C$ the Weyl chamber --
a fundamental domain for the shifted action attached to the fixed
triangular decomposition; $P^{+}=P\cap C$;

the $\co$ category of $\gtg$-modules attached to
the triangular decomposition;

a Verma module $M_{\lambda}\in\co,\; \lambda
\in\gth$ and a simple finite dimensional
module
 $V_{\lambda},\; \lambda\in P^{+}\subset P$;

the affine Lie algebra $\hgtg=\nc((z))\oplus\nc K$ and the ``generalized''
Borel subalgebra $\hgtg_{\geq}=\gtg\otimes\nc[[z]]\oplus\nc K$;

$\co_{k}$ -- the category of modules at level $k$
 (i.e. $K\mapsto k$),
and the full subcategory $\tco\subset\co_{k}$ consisting of
$\hgtg$-modules semisimple over $\gtg\subset\hgtg$;

$M_{\lambda}^{k}\in\co_{k},\; \lambda\in\gth^{\ast}$ is a Verma module;
 $V_{\lambda}^{k}=Ind_{\hgtg_{\geq}}^{\hgtg}\in\tco,\;\lambda\in P^{+}$ is a
Weyl module; more generally,
if $V$ is a $\gtg$-module, then $V^{k}\in\tco$ is a $\hgtg$-module
obtained by inducing from $V$; obviously, $V_{\lambda}^{k}$ is a
quotient of $M_{\lambda}^{k}$; each simple module is a quotient
of $M_{\lambda}^{k}$; denote it by $L_{\lambda}^{k}$;

if $k\not{\in}\nq$, then $\tco$ is semi-simple, each object being
a direct sum of Weyl modules; there is an obvious analogue of this
statement for $\co_{k}$;

for $k+h^{\vee}=p/q\in\nq_{\geq}$ consider an affine Weyl group
$W_{k}=pQ\propto W$, where $Q$ is a root lattice of $\gtg$; there
is the usual and the dotted (shifted) action of $W_{k}$ on
$\gth^{\ast}$; the fundamental domain for the latter is
$C_{aff}=C\cap\{\lambda: 0<(\lambda+\rho,\theta)<p\}$, where
$\theta$ is the highest root of $\gtg$; set
$P^{+}_{k}=P^{+}\cap C$; call $\lambda\in P^{+}_{k}$ (sometimes
$(\lambda,k)$ if $\lambda$ satisfies this condition) dominant;

by \cite{deodgabbkac,kac_kazhd,rcw}, $\tco=\oplus_{\lambda\in P^{+}_{k}}
\tco^{\lambda}$, where $\tco^{\lambda}$ is a full subcategory
consisting of modules whose composition series contain only
irreducible modules $L_{w\cdot\lambda}^{k}, w\in W_{k}$; similar
decomposition is true for $\co_{k}$.

\bigskip

{\em Duality Functors.} Given a vector space $W$, denote by $W^{d}$ its
total dual. If $W$ is a Lie algebra module, then so is $W^{d}$.

 Given a vector space $W$ carrying  a gradation by finite
dimensional subspaces, denote by $D(W)$ its restricted dual.

Objects of $\tco$ are canonically graded. Denote by $D:\co_{k}\rightarrow
\tco$, $M\mapsto D(M)$ the functor such that the $\hgtg$-module structure
is defined by precomposing the canonical action on the dual space with
an automorphism $\hgtg\rightarrow\hgtg$,
$g\otimes z^{n}\mapsto g\otimes(-z)^{-n}$.

The functors $^{d}, D(.)$ are exact.

There is an involution $\bar{ }: P^{+}\rightarrow P^{+}$ so that
$V_{\lambda}^{d}=V_{\bar{\lambda}}$.

\subsection{Geometry of weights}
The following is proved in \cite{jantz} Lemma 7.7.
\begin{lemma}
\label{lemmaongeomofw}
Suppose:

(i) $(\lambda,k),(\mu,k)\in P^{+}_{k}$ are regular;

(ii)$\bar{w}\in W$ satifies $\bar{w}(\lambda-\mu)\in P^{+}$;

(iii) $\nu$ is a weight of $V_{\bar{w}(\lambda-\mu)}$ such that
$w_{1}\cdot\lambda=w\cdot\mu +\nu$ for some $w,w_{1}\in W_{k}$.

Then: $w_{1}=w$ and $\nu\in W(\lambda-\mu)$.
\end{lemma}

\section{ the kazhdan-lusztig tensoring}
\label{kazhdluszttens}

Kazhdan and Lusztig \cite{kazh_luszt_0,kazh_luszt,kazh_luszt_1} (inspired by
Drinfeld \cite{drinf_1}) defined
a covariant bifunctor

\begin{equation}
\label{notattens}
\tco\times\tco\rightarrow\tco,\; A,B\mapsto A\tens B.
\end{equation}

We shall review its definition and main properties.
\subsection {Definition}
\label{twodefinitions}

\subsubsection{ The set-up}
\label{thesetup}

The notation to be used is as follows:

 $z$ is a  once and for all fixed coordinate on $\cp$;

$L\gtg^{P}, P\in\cp$ is the loop algebra attached to $P$; in other words,
$L\gtg^{P}=\gtg\otimes\nc((z-P)), P\in\nc$, and
$L\gtg^{\infty}=\gtg\otimes\nc((z^{-1}))$;

more generally, if $P=\{P_{1},...,P_{m}\}\subset\cp$, then
 \[L\gtg^{P}=\oplus_{i=1}^{m} L\gtg^{P_{i}};\]

$\hgtg^{P}=L\gtg^{P}\oplus \nc K, P\in\cp$ is the affine algebra attached to
the point $P$ --
the canonical central extension of $L\gtg^{P}$; of course, $\hgtg^0=\hgtg$;

more generally, if $P=\{P_{1},...,P_{m}\}\subset\cp$, then
 $\hgtg^{P}$ is the direct sum of  $\hgtg^{P_{i}},\; i=1,...,m$ modulo the
relation: all canonical central elements $K$ (one in each copy) are equal
each other;

$\Gamma=\gtg\otimes\nc[z,z^{-1},(z-1)^{-1}]$; $\Gamma$ is obviously a Lie
algebra.

The Laurent series expansions at points $\infty,1,0$ produce the Lie algebra
homomorphism
\[\epsilon: \Gamma\rightarrow L\gtg^{\{\infty,1,0\}}.\]

\begin{lemma}
\label{trehtochech}
The map $\epsilon$ lifts to a Lie algebra homomorphism
\[ \Gamma\rightarrow \hgtg^{\{\infty,1,0\}}.\]
\end{lemma}

Proof consists of using the residue theorem, see \cite{kazh_luszt}.

By pull-back, any $\hgtg^{\{\infty,1,0\}}$-module is canonically a
$\Gamma$-module.
 Further, any $A\in\tco$ is canonically a $\hgtg^{P}$-module
for any $P$  -- by the obvious change of coordinates; refer to this
as attaching $A$ to $P\in\cp$. Given $A,B,C\in\tco$, we shall regard $A\otimes
B\otimes C$
as a  $\hgtg^{\{\infty,1,0\}}$-module meaning that  $\hgtg^{\infty}$ acts on
$A$,  $\hgtg^{1}$ on $B$,  $\hgtg^{0}$ on $C$. (There is an obvious
ambiguity in this notation.) There arises the space of coinvariants
 \[(A\otimes B\otimes C)_{\Gamma}= (A\otimes B\otimes C)/\Gamma(A\otimes
B\otimes C).\]

This constrcuction easily generalizes to the case when instead of three points
-- $\infty,1,0$ -- there are $m$ points,   $m$ modules and instead of $\Gamma$
one considers the Lie algebra of rational functions on $\cp$ with $m$ punctures
with values in $\gtg$. We shall be mostly interested in the case $m=3$
and sometimes
in the case $m=2$. If $m=2$, then $\Gamma$ becomes
 $\tgtg=\gtg\otimes\nc[z,z^{-1}]$.

\begin{lemma}
\label{morphbetw2andcoinv}
Suppose $D(B)$ is attached to $\infty$, $A$ to 0. Then
\[Hom_{\hgtg}(A,B)=((D(B)\otimes A)_{\tgtg})^{d}.\]
\end{lemma}

Proof can be found in \cite{kazh_luszt}; the reader may also observe that
the arguments from \ref{morphismsandcoinvariants} are easily adjusted to this
case.

\subsubsection{Definition}
\label{Theseconddefinition}
Let $\hg$ be the central extension of $\Gamma$, the cocycle being defined
as usual except that one takes the sum of residues at $\infty$ and 1. Let
$\Gamma(0)\subset\hg$ be the subalgebra consisiting of functions vanishing
at 0. Obviously, $\Gamma(0)$ can also be regarded as a subalgebra of $\Gamma$.

Consider the (total) dual space  $(A\otimes B)^{d}$; it is naturally a
 $\hg$-module.  $(A\otimes B)^{d}$ carries the increasing filtration
$\{ (A\otimes B)^{d}(N)\}$, where

\begin{multline}
(A\otimes B)^{d}(N)\\
=\{x\in (A\otimes B)^{d}: \gamma_{1}\cdots\gamma_{N}x=0 \mbox{ if all }
\gamma_{i}\in \Gamma(0),
x\in (A\otimes B)\}.
\end{multline}

The space
$\cup_{N\geq 1} (A\otimes B)^{d}(N)$ is naturally a $\hgtg$-module. The passage
from $(A\otimes B)^{d}$ to $\cup_{N\geq 1} (A\otimes B)^{d}(N)$ (or its obvious
versions) is  often called  a functor of smooth vectors.

{\bf Define}
\begin{equation}
\label{seconddefform}
A\tens B=D(\bigcup_{N\geq 1} (A\otimes B)^{d}(N)).
\end{equation}

\begin{lemma}
\label{rightexactness}
The functor $\tens : \tco\times\tco\rightarrow \tco$ is right exact in
each variable.
\end{lemma}

{\bf Proof}  (see {\em loc. cit.}) The functor $\tens$ is a composition of two
dualizations, $^d$ and $D(.)$,
and the functor of smooth vectors. It is enough to remark that the first two
are exact while the last is only left exact.$\qed$

\subsection {Some properties of $\tens$.}

\subsubsection{ }
For the future reference we collect some of the properties
of $\tens$ in the following

\begin{theorem}
\label{mainthonklus}

(i)
\[Hom_{\hgtg}(A\tens B, D(C))= Hom_{\hgtg}(C,D(A\tens B))
=((A\otimes B\otimes C)_{\Gamma})^{d}.\]

(ii) If $A,B\in\tco$ have a Weyl filtration, then $A\tens B$ has
also. (Here by Weyl filtration we mean a filtration such that its quotients
are Weyl modules.)

(iii) If $k\not\in\nq$, then $V_{\lambda}^{k}\tens V_{\mu}^{k}
=(V_{\lambda}\otimes V_{\mu})^{k}$.

(iv) For any $k\in\nc$, $V_{\lambda}^{k}\tens V_{\mu}^{k}$ has
 a Weyl filtration (see(ii)), the multiplicity of $V_{\nu}^{k}$
being equal $(V_{\lambda}\otimes V_{\mu}: V_{\nu})$ (c.f. (iii)).

(v) There is an isomorphism $A\tens V_{0}^{k}\rightarrow A$
for any $A\in\tco$.

(vi) There are commutativity and asociativity morphisms
$A\tens B\approx B\tens A$ and $(A\tens B)\tens C\approx A\tens (B\tens C)$
which endow $\tco$ with the structure of a braided monoidal category.

\end{theorem}

\subsubsection{Morphisms and coinvariants.}
\label{morphismsandcoinvariants}
 The description
of morphisms in terms of coinvariants (see Theorem \ref{mainthonklus}(i))
is the hallmark of this theory. Let us briefly explain why (i) holds.
There is the obvious isomorphism of vector spaces
\[(A\otimes B\otimes C)^{d}\rightarrow  Hom_{\nc}(C,(A\otimes B)^{d}).\]
It induces the map

\[((A\otimes B\otimes C)_{\Gamma}^{d}\rightarrow  Hom_{\hg}(C,(A\otimes
B)^{d}).\]
By $\hg$-linearity, it actually gives the map

\[((A\otimes B\otimes C)_{\Gamma}^{d}\rightarrow  Hom_{\hg}(C,\bigcup_{N\geq 1}
(A\otimes B)^{d}(N)).\]
It remains to look at (\ref{seconddefform}) and note that $\hg$ is dense in
$\hgtg$.

\subsubsection{Using the spaces of coinvariants.}
\label{Usingthespacesofcoinvariants}
 A lot about the functor $\tens$ easily
follows from Theorem \ref{mainthonklus}(i). As an example, let us derive  (v).
By (i),

\[Hom_{\hgtg}(A\tens V_{0}^{k},B)=((A\otimes V_{0}^{k}\otimes
D(B))_{\Gamma})^{d}
\mbox{ for any } B\in\tco.\]

As $V_{0}^{k}=Ind_{\hgtg_{\geq}}^{\hgtg}\nc$, the Frobenius reciprocity
gives

\[(A\otimes V_{0}^{k}\otimes B)_{\Gamma}=(A\otimes D(B))_{\tgtg},\]
the latter space being  $Hom_{\hgtg}(A,B)$ by Lemma \ref{morphbetw2andcoinv}.
We see that the spaces
of morphisms of the modules $A$ and $A\tens V_{0}^{k}$ are equal, hence
so are the modules.

Replacing in this argument $\nc$ with a suitable finite dimensional
$\gtg$-module and repeating it three times one gets

\begin{equation}
\label{morphofweylintodual}
Hom_{\hgtg}(V_{\lambda}^{k}\tens V_{\mu}^{k}, D(V_{\nu}^{k}))
=Hom_{\gtg}(V_{\lambda}\otimes V_{\mu}, V_{\bar{\nu}}).
\end{equation}

As for generic $k$ $D(V_{\nu}^{k})\approx V_{\nu}^{k}$
(see \ref{prelim_reps}), (\ref{morphofweylintodual}) along with
Theorem \ref{mainthonklus}(i) implies Theorem \ref{mainthonklus}(iii).

\section{Affine translation functors}

\subsection{Definition.}
\label{defoftrfunctor}
For any  $(\lambda,k )\in P^{+}_{k}$ denote by $\tco^{\lambda}$ the full
subcategory
of $\tco$ consisting of modules whose composition factors all have highest
weights
 lying in the orbit $W_{k}\cdot(\lambda,k)$. There arises the projection
\[p_{\lambda}:\tco\rightarrow \tco^{\lambda}.\]
This all has been reviewed in \ref{prelim_reps}.

Given  $(\lambda,k ),(\mu,k)\in P^{+}_{k}$, pick $\bar{w}\in W$ so that
$\bar{w}(\lambda-\mu)\in P^{+}$. It is easy to see that then
$(\bar{w}(\lambda-\mu),k)\in P^{+}_{k}$.

{\bf Define } the translation functor

\begin{equation}
\label{defoftrfunctor_form}
T_{\mu}^{\lambda}:\; \tco^{\mu}\rightarrow \tco^{\lambda}
A\mapsto p_{\lambda}(V_{\bar{w}(\lambda-\mu)}^{k}\tens A).
\end{equation}
This functor was first introduced by Finkelberg \cite{fink} who, however,
considered it only
for $k<0$.

As an immediate corollary of the definition, one has
\begin{equation}
\label{howtoadjointtlm}
T_{\lambda}^{\mu}=p_{\mu}\circ ((V_{\bar{w}(\lambda-\mu)}^{d})^{k}\tens ?)
\end{equation}

\subsection{Rigidity of  Weyl modules with dominant highest weight.}
\label{rigidityofweylmod}
\begin{lemma}
\label{actoftlmonweyls}

 If $(\lambda,k),(\mu,k)$ are regular (i.e. off the affine walls) and $w\in
W_{k}$
satisfies $w\cdot\mu\in P^{+}$, then
\[T_{\mu}^{\lambda}(V_{w\cdot\mu}^{k})=V_{w\cdot\lambda}^{k}.\]

\end{lemma}

{\bf Proof.}  By Theorem \ref{mainthonklus} (iv),
$T_{\mu}^{\lambda}(V_{w\cdot\mu}^{k})$
has a filtration with quotients isomorphic to $V_{w_{1}\cdot\lambda}^{k}$,
$w_{1}\in W_{k}$ such that $w_{1}\cdot\lambda =w\cdot\mu+\nu$, $\nu$ being
a weight of $V_{\bar{w}(\lambda-\mu)}$. By Lemma \ref{lemmaongeomofw},
$w_{1}=w$. This implies that this filtration has only one term,
$V_{w\cdot\lambda}^{k}$. $\qed$

\begin{corollary}
\label{rigofvlk}
If $(\lambda,k)\in P^{+}_{k}$ is regular, then $V_{0}^{k}$ is a direct summand
of
$V_{\lambda}^{k}\tens V_{\bar{\lambda}}^{k}$.
\end{corollary}

{\bf Proof.} Of course $(0,k)$ is dominant regular
and $p_{0}A$ is a direct summand of $A$. It remains
to observe that $T_{\bar{\lambda}}^{0}V_{\bar{\lambda}}^{k}=
p_{0}(V_{\lambda}^{k}\tens V_{\bar{\lambda}}^{k})$ and use
 Lemma \ref{actoftlmonweyls} to get
 $T_{\bar{\lambda}}^{0}V_{\bar{\lambda}}^{k}=
V_{0}^{k}$ . $\qed$

We get the maps
\[i_{\lambda}: V_{0}^{k}\rightarrow V_{\lambda}^{k}\tens
 V_{\bar{\lambda}}^{k},\;
e_{\lambda}: V_{\bar{\lambda}}^{k}\tens V_{\lambda}^{k}\rightarrow V_{0}^{k}.\]

Observing that the maps between $\tens$-products of Weyl
modules are uniquely determined by the induced maps of
the corresponding finite dimensional $\gtg$-modules (Theorem \ref{mainthonklus}
and (\ref{morphofweylintodual}) ), we see that
 we can normalize $i_{\lambda},e_{\lambda}$ so that the compositions

\begin{multline}
V_{\lambda}^{k}=V_{0}^{k}\tens V_{\lambda}^{k}\stackrel{i_{\lambda}\otimes
id}{\rightarrow}
V_{\lambda}^{k}\tens V_{\bar{\lambda}}^{k}\tens V_{\lambda}^{k}
\stackrel{id\otimes e_{\lambda}}{\rightarrow} V_{\lambda}^{k}\\
V_{\bar{\lambda}}^{k}=V_{\bar{\lambda}}^{k}\tens V_{0}^{k}\stackrel
{id\otimes i_{\lambda}}
{\rightarrow}
V_{\bar{\lambda}}^{k}\tens V_{\lambda}^{k}\tens V_{\bar{\lambda}}^{k}
\stackrel{ e_{\lambda}\otimes id}{\rightarrow} V_{\bar{\lambda}}^{k},
\label{feofrigid}
\end{multline}

are equal to the identity. By definition (see e.g. \cite{kazh_luszt_1} III,
Appendix) we have

\begin{corollary}
\label{rigidity}
If $(\lambda, k)\in P^{+}_{k}$, then $V_{\lambda}^{k}$ and
$V_{\bar{\lambda}}^{k}$
are rigid.
\end{corollary}

Consider the functor $V_{\lambda}^{k}\tens ?: \tco\rightarrow\tco,
 M\mapsto V_{\lambda}^{k}
\tens M.$

\begin{corollary}
\label{adjointnessproperty}
(i)If $(\lambda, k)\in P^{+}_{k}$, then the functors $V_{\lambda}^{k}\tens ?$
 and $V_{\bar{\lambda}}^{k}\tens ?$ are adjoint, i.e. there is a functor
ismorphism

\[Hom_{\hgtg}(V_{\lambda}^{k}\tens A,B)=Hom_{\hgtg}(
A,V_{\bar{\lambda}}^{k}\tens B).\]

(ii) If $(\lambda, k)\in P^{+}_{k}$, then the functors $V_{\lambda}^{k}\tens ?$
 and $V_{\bar{\lambda}}^{k}\tens ?$ are exact, i.e. send exact short sequences
to
exact ones.
\end{corollary}

{\bf Proof} is standard; for the reader's convenience we reproduce the one from
 \cite{kazh_luszt_1} III,
Appendix.
To prove (i), consider two composition maps

\[\phi:Hom_{\hgtg}(V_{\lambda}^{k}\tens A,B)\rightarrow
Hom_{\hgtg}(V_{\bar{\lambda}}^{k}\tens  V_{\lambda}^{k}\tens
A,V_{\bar{\lambda}}^{k}\tens B)
\stackrel{i_{\bar{\lambda}}}{\rightarrow}
Hom_{\hgtg}( A,V_{\bar{\lambda}}^{k}\tens B),\]

\[\psi:Hom_{\hgtg}( A,V_{\bar{\lambda}}^{k}\tens B)\rightarrow
Hom_{\hgtg}(V_{\lambda}^{k}\tens A,V_{\lambda}^{k}\tens
V_{\bar{\lambda}}^{k}\tens B)
\stackrel{e_{\bar{\lambda}}}{\rightarrow }
Hom_{\hgtg}(V_{\lambda}^{k}\tens A, B).\]

By (\ref{feofrigid}), the compositions $\phi\circ\psi$ and $\psi\circ\phi$ are
equal to
the identity.

(ii) is an easy consequence of (i): we have to prove that $B_{1}\rightarrow
B_{2}$ is
a monomorphism implies that $V_{\lambda}^{k}\tens B_{1}\rightarrow
V_{\lambda}^{k}\tens B_{2}$ is also, or, equivalently, that for any $A\in\tco$
the induced map
\[Hom_{\hgtg}(A,V_{\lambda}^{k}\tens B_{1})\rightarrow
Hom_{\hgtg}(A,V_{\lambda}^{k}\tens  B_{2})\]
is also a monomorphism. By (i), it is equivalent to proving that
\[Hom_{\hgtg}(V_{\bar{\lambda}}^{k}\tens A, B_{1})\rightarrow
Hom_{\hgtg}(V_{\bar{\lambda}}^{k}\tens A, B_{2})\]
is a monomorphism, but this is an obvious corollary of injectivity of the map
$B_{1}\rightarrow B_{2}$. $\qed$

\subsection{Properties of affine translation functors.}
Recall that there is the notion of a formal character $ch A$ for
any $A\in\tco^{\lambda}$, see e.g. \cite{deodgabbkac}. There arises an abelian
group of characters, each of the following sets being a topological basis in
it:

$\{ch V_{w\cdot\lambda}^{k},\; w\in W_{k}\}$,
$\{ch L_{w\cdot\lambda}^{k},\; w\in W_{k}\}$.
Of course the symbols $ch V_{w\cdot\lambda}^{k},\; ch L_{w\cdot\lambda}^{k}$
should be ignored unless $w\cdot\lambda\in P^{+}$. Observe that
\begin{equation}
\label{mezhduvil}
ch A= \sum_{w\geq w_{o}}\bar{n}_{w}ch L_{w\cdot\mu}^{k}\Leftrightarrow
ch A= \sum_{w\geq w_{o}}n_{w}ch V_{w\cdot\mu}^{k}
\end{equation}

\begin{proposition}
\label{conseqreltochar}
Let  $(\lambda,k),(\mu,k)$ be regular dominant.

(i) $T_{\mu}^{\lambda}$ is exact;

(ii)  $T_{\mu}^{\lambda}, T_{\lambda}^{\mu}$ are adjoint to each other;

(iii)  If $ch A=\sum_{w\in W_{k}}n_{w}ch V_{w\cdot\mu}^{k}$, then
 $ch T_{\mu}^{\lambda} A=\sum_{w\in W_{k}}n_{w}ch V_{w\cdot\lambda}^{k}$.

(iv) $T_{\mu}^{\lambda}(L_{w\cdot\mu}^{k})=L_{w\cdot\lambda}^{k}$;

(v) More generally,  $T_{\mu}^{\lambda}(.)$ establishes
an equivalence of the submodule lattices of
$V_{w\cdot\mu}^{k}$ and $V_{w\cdot\lambda}^{k}$.
\end{proposition}

{\bf Proof.}
(i)  $T_{\mu}^{\lambda}$ is exact as a composition of the exact functors
$p_{\lambda}$ and $V_{\bar{w}(\lambda-\mu)}^{d})^{k}\tens ?$, see
Corollary \ref{adjointnessproperty} (ii).

(ii) By Corollary \ref{adjointnessproperty} (i),
one has for any $A\in\tco^{\mu},\; B\in\tco^{\lambda}$

\begin{multline}
Hom_{\hgtg}(T_{\mu}^{\lambda}A,B)=
Hom_{\hgtg}(p_{\lambda}(V_{\bar{w}(\lambda-\mu)}^{k}\tens A),B)
=Hom_{\hgtg}(V_{\bar{w}(\lambda-\mu)}^{k}\tens A,B)\\=
Hom_{\hgtg}( A,(V_{\bar{w}(\lambda-\mu)}^{d})^{k}\tens B)
=Hom_{\hgtg}( A,p_{\mu}(V_{\bar{w}(\lambda-\mu)}^{d})^{k}\tens B)\\
=Hom_{\hgtg}( A, T_{\lambda}^{\mu}B).
\end{multline}

(iii) follows at once from (i) ( if one uses the local composition series, see
e.g. \cite{deodgabbkac}).

(iv) Let $T_{\mu}^{\lambda}(L_{w_{0}\cdot\mu}^{k})$ be reducible. There
arises an exact sequence with non-zero $N$
\[0\rightarrow N\rightarrow  T_{\mu}^{\lambda}(L_{w_{0}\cdot\mu}^{k})
\rightarrow L_{w_{0}\cdot\lambda}^{k}\rightarrow 0.\]

Applying $T_{\lambda}^{\mu}$ to it one gets

\[0\rightarrow T_{\lambda}^{\mu}( N)\rightarrow
T_{\lambda}^{\mu}(T_{\mu}^{\lambda}(L_{w_{0}\cdot\mu}^{k}))
\rightarrow T_{\lambda}^{\mu}(L_{w_{0}\cdot\lambda}^{k})\rightarrow 0.\]

By (iii) and (\ref{mezhduvil}),
 $ch(T_{\lambda}^{\mu}(T_{\mu}^{\lambda}(L_{w_{0}\cdot\mu}^{k})))=
ch L_{w_{0}\cdot\mu}^{k}$ and $ch T_{\lambda}^{\mu}( N)\neq 0$; therefore
$ch T_{\lambda}^{\mu}(L_{w_{0}\cdot\lambda}^{k})< ch L_{w_{0}\cdot\mu}^{k}$.
Contradiction.

(v) Here proof is an obvious version of that of (iv). By
  (ii) it is enough to show that
if $A\subset B\subset V_{w\cdot\mu}^{k}$, then
 $T_{\mu}^{\lambda}(A)\subset T_{\mu}^{\lambda}( B)\subset V_{w\lambda}^{k}$.
Using (\ref{mezhduvil}) and passing to quotients, if necessary, the
problem is reduced to the case when $B$ is a highest weight module. In this
case the arguments of (ii) go through practically unchanged.
$\qed$

\subsection { }

\begin{theorem}
\label{equivofcat}
The functor $T_{\mu}^{\lambda}:\tco^{\mu}\rightarrow\tco^{\lambda}$ is
an equivalence of categories.
\end{theorem}

{\bf Proof.} It is enough show that $T_{\mu}^{\lambda}\circ T_{\lambda}^{\mu}:
\tco^{\lambda}\rightarrow \tco^{\lambda}$ is equivalent to the identity.
In other words, we want to show that $id: A\rightarrow A, \;
A\in\tco^{\lambda}$ is
transformed into an isomorphism in
 $Hom_{\hgtg}(T_{\mu}^{\lambda}\circ T_{\lambda}^{\mu}(A), A)$. We already know
this when $A$
 is simple, see
Corollary \ref{conseqreltochar} (ii). Using the local composition
series one proves it for an arbitrary $A$.

An alternative way to prove the theorem is to observe that by Corollary
\ref{rigofvlk} the action of  $T_{\mu}^{\lambda}\circ T_{\lambda}^{\mu}$ is
equivalent to that of $V_{0}^{k}\tens ?$, the latter being equivalent
to $id$ by Theorem \ref{mainthonklus} (v). $\qed$

\subsection { Generalizing from $\tco$ to $\co_{k}$.}
\label{genfromtcotook}

Our two key resluts -- Proposition \ref{conseqreltochar} and Theorem \ref
{equivofcat} -- can be carried over to the category $\co_{k}$. Let us briefly
explain it. We will be using subcategories $\co_{k}^{\lambda}\subset
\co_{k}$ (see \ref{prelim_reps}) only when
$k+h^{\vee}\in\nq_{>}$ and $\lambda$ is integral, although the last condition
can be easily relaxed.

It is non-trivial (if at all meaningful ) to carry the Kazhdan-Lusztig
tensoring over
 to the entire $\co_{k}$. (An intelligent way to do something like it
requires introducing additional structures, see \cite{feimal2}.) It is
however straightforward to extend it to the functor
\[\tens:\;\tco\times\co_{k}\rightarrow \co_{k},\]
as proposed by Finkelberg \cite{fink}.  One basic property of this operation
absolutely analogous (along with the proof) to Theorem \ref{mainthonklus} (iv)
is as follows.

As $V_{\lambda}\otimes M_{\mu}$ has a filtration by Verma modules in the
category of
$\gtg$-modules, $V_{\lambda}^{k}\tens M_{\mu}^{k}$ has a filtration
by Verma modules in $\co_{k}$; further the multipliciites are the same as
in the finite dimensional case:
\[(V_{\lambda}^{k}\tens M_{\mu}^{k}:M_{\nu}^{k})=
(V_{\lambda}\otimes M_{\mu}: M_{\nu}).\]

Given this one can easily inspect our exposition of affine translation
functors and
observe that quite a lot  carries over to the setting of
$\co_{k}$ word for word except that at the
appropriate places Weyl modules are to be changed for the corresponding
Verma modules. Here are some examples:

(i) definition of $T_{\mu}^{\lambda}:\co_{k}^{\mu}\rightarrow\co_{k}^{\lambda}$
if $\lambda,\;\mu$  belong to the same Weyl chamber;

(ii) the Verma filtration of $V_{\lambda}^{k}\tens M_{w\cdot\mu}^{k}, w\in
W_{k}$ and
 Lemma \ref{lemmaongeomofw} imply that $T_{\mu}^{\lambda}(M_{w\cdot\mu}^{k})=
M_{w\cdot\lambda}^{k}$ if $(\mu,k),(\lambda,k)$ are regular
 (c.f. Lemma \ref{actoftlmonweyls}); observe that we can now drop the condition
that
$w\cdot\mu\in P^{+}$;

(iii) therefore Proposition \ref{conseqreltochar} holds with the indicated
changes.

 We get
\begin{theorem}
\label{eqofcat_forbigg}
The functor $T_{\mu}^{\lambda}:\co^{\mu}_{k}\rightarrow\co^{\lambda}_{k}$ is
an equivalence of categories if
 $\lambda,\;\mu$ are integral and both belong to the same Weyl chamber.
\end{theorem}

\section {annihilating ideals of highest weight modules}

\subsection{Vertex operators and ...}

The usual tensor functor $\otimes: M,N\mapsto M\otimes N$ has
the following fundamental(and trivial) property: there is a natural map
\begin{multline}
\label{fundmapatttoustens}
N\rightarrow Hom_{\nc}(M,M\otimes N)\\
n\mapsto n(.) \mbox{ such that } n(m)=m\otimes n.
\end{multline}

Here  we shall explain the $\tens$-analogue of this map

\subsubsection{ }
\label{Vertexoperatorsandvertexoperatoralgebras}

By Theorem \ref{mainthonklus} (v), $A\tens V_{0}^{k}\approx A$ for any
$A\in\tco$.
Therefore by Theorem \ref{mainthonklus} (i), there is a natural isomorphism

\[
((A\otimes V_{0}^{k}\otimes D(A))_{\Gamma})^{d}\approx Hom_{\hgtg}(A,B), \]
for any $B\in\tco$.

 Recall that the  space
$((A\otimes V_{0}^{k}\otimes D(A))_{\Gamma})^{d}$ was defined by means of
$\Gamma$, the
latter being defined by choosing three points,  $ \infty,1,0$,
see the end of \ref{Usingthespacesofcoinvariants}. The
choice of points was, of course, rather arbitrary. Keeping $ \infty, 0$
fixed and $A$, $D(B)$ attached to $\infty,0$ resp., we shall allow the third
point to vary. We get then the family
of Lie algebras $\Gamma_{t}, t\in\nc^{\ast}$ and the family
of the  one-dimensional spaces (c.f. \ref{thesetup})
\[<A,V_{0}^{k},D(B)>_{t}:=
((A\times V_{0}^{k}\times D(B))_{\Gamma_{t}})^{d},\; t\in\nc^{\ast}. \]

These naturally arrange in a trivial line bundle over $\nc^{\ast}$,
the fiber being isomorphic to
\[<A,V_{0}^{k},D(B)>_{t}=(A\otimes D(B))_{\tgtg}=
Hom_{\hgtg}(A,B),\]
by the arguments using  Frobenius reciprocity as in
\ref{Usingthespacesofcoinvariants}. Pick a section of this bundle by choosing
$\phi\in Hom_{\hgtg}(A,B)$.

 Hence we get a trilinear  functional (depending on $t\in\nc^{\ast}$)
\[\Phi_{t}^{\phi}\in <A,V_{0}^{k},D(B)>_{t}\subset (A\otimes V_{0}^{k}\otimes
D(B))^{d}.\]

Reinterprete it as the linear map:

\begin{equation}
\label{fromcoinvtomaps_0}
\Phi_{t}^{\phi}(.): V_{0}^{k}\rightarrow (A\otimes  D(B))^{d},
\end{equation}

or, equivalently,
\begin{equation}
\label{fromcoinvtomaps}
\tilde{\Phi}_{t}^{\phi}(.): V_{0}^{k}\rightarrow Hom_{\nc}(A,  D(B)^{d}),\;
t\in\nc^{\ast}.
\end{equation}

The latter map is an analogue of $N\rightarrow Hom_{\nc}(M,M\otimes N)$
 mentioned above.
To analyze its properties observe that there is an obvious embedding
$B\rightarrow (D(B))^{d}$. It does not, of course,
 allow us to interprete $\tilde{\Phi}_{t}^{\phi}(v),\;v\in V_{0}^{k}$
as an element of $Hom_{\nc}(A,B)$ depending on $t$. But, as the following
lemma shows, Fourier coefficients of
 $\tilde{\Phi}_{t}^{\phi}(v),\;v\in V_{0}^{k}$ are actually elements
of $Hom_{\nc}(A,B)$.
To formulate this lemma observe that there is a natural gradation
on $A$ and $B$ consistent with that of $\tgtg$; e.g.
$A=\oplus_{n\geq 0}A[n],\; dim A[n]<\infty$.

\begin{lemma}
\label{whycantfour}

Let $B$ be either $A$ or a quotient of $A$, $id: A\rightarrow B$ be the natural
projection. Then:

(i) $\Phi_{t}^{id}(vac)(x,y)=y(x),$ where $vac$ is understood as the generator
of $V_{0}^{k}$;

(ii) more generally, if $v\in V_{0}^{k}[n],\; x\in A[m],\;y\in D(B)[l]$,
then
\[\Phi_{t}^{id}(v)(x,y)\in \nc\cdot t^{-l+m-n}.\]
\end{lemma}

{\bf Proof.}  Given $g\in\gtg$, denote by $g_{n}\in\hgtg^{P}$ the
element $g\otimes (z-P)^{n}$ or $g\otimes z^{-n}$ if
$P=\infty$. (It should be clear from the
context which $P$ is meant.) Thus $g_{n}x=(g\otimes z^{-n})x$ if $x\in A$, the
$A$
 being attached to $\infty$; similarly,
  $g_{n}x=(g\otimes z^{n})x$ if $x\in D(B)$, the $D(B)$ being attached to $0$.

 (i) can be proved by an obvious induction on the
degree of $x$ and $y$ using the following formula (which follows from the
definition
of  $((A\times V_{0}^{k}\times D(B))_{\Gamma})^{d}$ and the Laurent
expansions of $z^{-n}$ at $\infty$ and 0):

\[\Phi_{t}^{id}(vac)(g_{n}x,y)=-\Phi_{t}^{id}(vac)(x,g_{-n}y).\]

 To prove (ii) observe, first, that (i) is a particular case of (ii) when
$v=vac$. One then proceeds by induction on $n$ using the formula
 (which again follows from the definition
of  $((A\times V_{0}^{k}\times D(B))_{\Gamma})^{d}$
 and the Laurent
expansions of $(z-t)^{-n}$ at $\infty$ and 0):

\begin{multline}
(-1)^{n-1}(n-1)!\Phi_{t}^{id}(g_{-n}v)(x,y)\\
=(\frac{d}{dt})^{n-1}\{\sum_{i=1}^{\infty}t^{i-1}\Phi_{t}^{id}(v)(g_{i}x,y)
-\sum_{i=0}^{\infty}t^{-i-1}\Phi_{t}^{id}(v)(x,g_{i}y)\}.
\end{multline}

$\qed$

Observe that the spaces $A,B$ being graded, the space $Hom_{\nc}(A,D(B))$ is
also.
Lemma \ref{whycantfour} means that although the map
$\tilde{\Phi}_{t}^{id}(.)$ from (\ref{fromcoinvtomaps}) cannot be interpreted
as an element of  $Hom_{\nc}(A,B)$, its Fourier components can  because they
are homogeneous. To compare with \cite{fren_lep_meur} introduce the following
notation: for any $v\in V_{0}^{k}[n]$ set

\begin{equation}
Y(v,t)=\sum_{i\in\nz}v_{i} t^{-i-n},
\label{defofyvt}
\end{equation}
where
\begin{equation}
v_{i}:=\oint \tilde{\Phi}_{t}^{id}(v)t^{i+n-1}\, dt:
A[l]\rightarrow B[l+i],
\label{defofyvtfourcomp}
\end{equation}
for all $l\geq 0$, and call the generating functions $Y(v,t)$ {\em fields}.
 For example, it easily follows from the formulae above
that
\begin{equation}
\label{current}
x(t):=Y(x_{-1}vac,t)=\sum_{i\in\nz}x_{i}t^{-i-1},
\end{equation}
producing the famous {\em current} $x(t)$. Another fact easily reconstructed
from the formulae above (especially from the proof of Lemma \ref{whycantfour})
is that

\begin{equation}
\label{succcurrent}
(-1)^{n-1}(n-1)!Y(x_{-n}v,t)= :x(t)^{(n-1)}Y(v,t):,
\end{equation}

where we set
\[  :x(t)^{(n-1)}Y(v,z):\;=(x(z)^{(n-1)})_{-}Y(v,t) +
Y(v,t)(x(z)^{(n-1)})_{+},\]
$(x(z)^{(n-1)})_{\pm}$ being defined as usual (see e.g. \cite{frzhu}).  It
follows that all fields are infinite
combinations of elements of $\hgtg$.

The expressions  $Y(v,t)$ are not only formal generating functions.
In this notation Lemma \ref{whycantfour} can be rewritten as follows.

\begin{corollary}
\label{actionoffields3points}
Under the assumptions of Lemma \ref{whycantfour},
\[\Phi_{t}^{id}(v)(x,y) = y(Y(v,t)x).\]
\end{corollary}

\subsubsection{   }
\label{twogeneralizations} The considerations of
\ref{Vertexoperatorsandvertexoperatoralgebras}
 are easily generalized as follows.
(We shall skip the proofs as they  essentially  repeat those in
\ref{Vertexoperatorsandvertexoperatoralgebras}.)

 Replace $V_{0}^{k}$ with $V_{\lambda}^{k}$ and pick $A,B\in\tco$
so that the space $<A,V_{\lambda}^{k},D(B)>_{t}\neq 0$. For any
$\phi\in <A,V_{\lambda}^{k},D(B)>_{t}$ we get a map

\begin{multline}
Y(.,t):\; V_{\lambda}^{k}\rightarrow Hom_{\nc}(A,B((t,t^{-1})))\\
V_{\lambda}^{k}\ni v\mapsto Y(v,t)=\sum_{i\in\nz}v_{i}t^{-i-\tilde{v}},\;
v_{i}\in  Hom_{\nc}(A,B).
\end{multline}

$Y(v,t),\; v\in V_{\lambda}^{k}$ is a generating function having all properties
its counterpart from \ref{Vertexoperatorsandvertexoperatoralgebras}
with one notable exception. Consider the ``upper floor'' of $V_{\lambda}^{k}$:
$V_{\lambda}\subset V_{\lambda}^{k}$. The Fourier components
of the fields $Y(v,t),\; v\in V_{\lambda},\; \lambda\neq 0$ generate a
$\hgtg$-submodule of
$ Hom_{\nc}(A,B)$ isomorphic to the loop module
$L(V_{\lambda})=V_{\lambda}\otimes \nc[z,z^{-1}]$. Strange as it may seem to
be,
if $\lambda=0$, then instead of $\nc[z,z^{-1}]$ this construction gives
simply $\nc$ -- this was explained above.

The embedding $L(V_{\lambda})\subset Hom_{\nc}(A,B)$ is called {\em a vertex
operator}.
It is easy to see that all vertex operators are obtained via the described
construction.

\subsection{...and vertex operator algebras}
\label{and_vertex_operator_algebras}

We now recall that a vertex operator algebra (VOA) is defined to be a graded
vector space $\bigcup_{i\in\nz} V[i],\; dim V_{i}<\infty$ along with a map
\[Y(.,t): V\rightarrow End(V)((t,t^{-1})),\]
satisfying certain axioms among which we mention {\em associativity}
and {\em commutativity} axioms, see e.g. \cite{fren_lep_meur,frzhu}. Similarly
one defines the notion of a module (submodule) over a VOA. A VOA is a module
over itself;
call an ideal of a VOA a submodule of a VOA as a module over itself. Observe
that it follows from the associativity axiom that the Fourier components of
fields $Y(v,t),\; v\in V$ close in a Lie algebra, $Lie(V)$. In this way,
an ideal of a VOA $V$ produces an ideal of $Lie (V)$ in the Lie algebra sense.
Not any ideal of $Lie(V)$ can be obtained in this way. Refer to such an ideal
{\em an ideal of }  $Lie(V)$ as {\em VOA}.

It follows from \cite{frzhu} that the constructions of
\ref{Vertexoperatorsandvertexoperatoralgebras} give:
$(V_{0}^{k}, Y(.,t))$ is a vertex
operator algebra and each $A\in\tco$ is a module over it. $Lie(V_{0}^{k})$
is habitually denoted $U(\hgtg)_{loc}$ and called {\em a local completion}
of $U(\hgtg)$, even though it is not an associative algebra!
A moment's thought shows that the ideal lattice of
$U(\hgtg)_{loc}$ as VOA is isomorphic with the submodule lattice of
$V_{0}^{k}$
as a $\hgtg$-module.

\subsection{ }

Here we prove the following theorem -- one of the main results of
 this paper.

\begin{theorem}
\label{theoronannideals}
Let $k\in\nq_{>}$, $(\lambda,k),(0,k)\in P^{+}_{k}$ be regular.
Denote by $\Omega(V_{\lambda}^{k})$ the submodule lattice of
 $V_{\lambda}^{k}$,
and by $\Omega(U(\hgtg)_{loc})$ the ideal lattice of $U(\hgtg)_{loc}$ as VOA at
the level $k$.
There is a lattice equivalence
\begin{multline}
\omega: \Omega(U(\hgtg)_{loc})\rightarrow \Omega(V_{\lambda}^{k}),\\
 \Omega(U(\hgtg)_{loc})\ni I\mapsto IV_{\lambda}^{k}.
\end{multline}
\end{theorem}

{\bf Proof.}

 First of all, by definition \ref{and_vertex_operator_algebras}  $\omega$ is
equivalently reinterpreted as a map of the submodule lattices
of the $\hgtg$-modules: $\omega:\; \Omega(V_{0}^{k})\rightarrow
\Omega(V_{\lambda}^{k})$.
In what follows we shall make use of this reinterpretation.

  Consider the translation functor:
$T_{0}^{\lambda}$. If $N\subset V_{0}^{k}$ is a submodule,
then on the one hand we have
\[T_{0}^{\lambda}(V_{0}^{k})=V_{\lambda}^{k}\tens V_{0}^{k}
(=V_{\lambda}^{k}),\]
and therefore
\[T_{0}^{\lambda}(V_{0}^{k}/N)=V_{\lambda}^{k}\tens ( V_{0}^{k}/N).\]

By Theorem \ref{mainthonklus} and Corollary \ref{actionoffields3points},
\begin{multline}
Hom_{\hgtg}(T_{0}^{\lambda}(V_{0}^{k}/N),?)=
<V_{0}^{k}/N,V_{\lambda}^{k},D(?)>_{t}\\
= Hom_{\hgtg}(V_{\lambda}^{k}/\omega(N),?).
\end{multline}

On the other hand, by Proposition \ref{conseqreltochar} (i)
 \[T_{0}^{\lambda}(V_{0}^{k}/N)= V_{\lambda}/T_{0}^{\lambda}(N).\]

We conclude immediately that $\omega(N)= T_{0}^{\lambda}(N)$.
It remains to recollect that $T_{0}^{\lambda}$ is an isomorphism of
the submodule lattices by Proposition \ref{conseqreltochar} (v).
$\qed$

\bigskip

An application of this result to annihilating ideals of
{\em admissible representations} is as follows. Recall that if
$k+h^{\vee}\in\nq_{>}$, $(\lambda,k)\in P^{+}_{k}$ is regular, then
$L_{\lambda}^{k}$ is called admissible \cite{kac_wak}.
$L_{\lambda}^{k}$ is an irreducible quotient of
$V_{\lambda}^{k}$ be a submodule $N_{\lambda}^{k}$
 generated by one singular vector,
see also \cite{kac_wak}.  By Theorem \ref{theoronannideals},
$\omega(N_{0}^{k})=N_{\lambda}^{k}$. We get

\begin{corollary}
\label{annidadmrepr}
The annihilating ideal of an admissible representation  equals
$Lie(N_{0}^{k})$; in particular, it is generated (as VOA) by one
singular vector of $V_{0}^{k}$.
\end{corollary}

{\em Remarks.}

(i) In the case $\gtg=\gtsl_{2}$, Corollary \ref{annidadmrepr}
follows from the more general results of \cite{feimal1}, see also
\cite{feimal2}.

(ii) If the Feigin-Frenkel conjecture on the singular support of
$L_{0}^{k}$ (theorem in the $sl_{2}$-case, see \cite{feimal2})
were correct, then Corollary \ref{annidadmrepr}
 would imply its validity for any admissible representation from
$\tco$ and thus would give a new example of rational conformal field
theory.

(iii) Another way to think of Corolary \ref{annidadmrepr} is
that $L_{0}^{k}$ is a VOA and $L_{\lambda}^{k}$ is a module over
it; in the $sl_{2}$-case, this point of view is adopted in
 \cite{adam,donglimas}.

\section{What is a Harish-Chandra bimodule over an affine
Lie algebra?}

\label{whcafm}

\subsection{ Restricted Harish-Chandra category.}
Our approach to defining affine Harish-Chandra bimodules
will heavily rely on the properties of affine translation
functors. We begin in the framework of the category $\tco$,
 see Proposition \ref{conseqreltochar} and Theorem
\ref{equivofcat}. Call a triple of weights $\lambda_{l},
\lambda_{r},\lambda\in P^{+}_{k}$ a translation datum if
$\lambda_{l}-\lambda\in W\cdot\lambda_{r}$.
There arises
the translation functor
$T_{\lambda}^{\lambda_{l}}=p_{\lambda_{l}}\circ(V_{\lambda_{r}}^{k}
\tens ?)$.

Let $Funct(\tco^{\lambda_{r}},\tco^{\lambda_{l}})$ be
 the category of functors
from $\tco^{\lambda_{l}}$ to $\tco^{\lambda_{r}}$.
There is a functor
\[\Phi:\tco\rightarrow Funct(\tco^{\lambda_{r}},\tco^{\lambda_{l}}),\]
\[ \Phi(A): B\mapsto p_{\lambda_{l}}(A\tens B).\]

Setting for the sake of breavity $\tilde{\cf}_{\lambda_{r}}^{\lambda_{l}}=
Funct(\tco^{\lambda_{r}},\tco^{\lambda_{l}})$, we get
the natural map
\[i:\; Hom_{\tilde{\cf}_{\lambda_{r}}^{\lambda_{l}}}(F,G)
\rightarrow Hom_{\hgtg}(F(V_{\lambda_{r}}^{k}),G(V_{\lambda_{r}}^{k})),\]
where $i(\psi)$ is simply the value of the functor morphism $\psi$
on $V_{\lambda_{r}}^{k}$.

\begin{conjecture}
\label{conjhowdefhch}
If $(\lambda_{r},\mu,\lambda_{l})$ and
 $(\lambda_{r},\nu,\lambda_{l})$ are translation data and
$A\in\tco^{\mu},\;B\in\tco^{\nu}$, then the map
\[Hom_{\tilde{\cf}_{\lambda_{r}}^{\lambda_{l}}}(\Phi(A),\Phi(B))
\rightarrow Hom_{\hgtg}(\Phi(A)(V_{\lambda_{r}}^{k}),
\Phi(B)(V_{\lambda_{r}}^{k})).
\]
is an isomorphism (c.f. Theorem 3.5 in \cite{bernst_gelf}).
\end{conjecture}

To provide a supporting evidence,
we   prove surjectivity in the case $\mu=\nu$. As
$(\lambda_{r},\mu,\lambda_{l})$ and
 $(\lambda_{r},\nu,\lambda_{l})$ are translation data,
$\Phi(A)(V_{\lambda_{r}}^{k})=T_{\mu}^{\lambda_{l}}(A)$ and
$\Phi(B)(V_{\lambda_{r}}^{k})=T_{\mu}^{\lambda_{l}}(B)$.
By Theorem
\ref{equivofcat} we get an isomorphism
\[T_{\lambda_{l}}^{\mu}: \; Hom_{\hgtg}(\Phi(A)(V_{\lambda_{r}}^{k}),
\Phi(B)(V_{\lambda_{r}}^{k}))\approx
Hom_{\hgtg}(A,B).\]

\begin{sloppypar}
It follows that any $\phi\in Hom_{\hgtg}(A,B)$
gives rise to $\Phi(\phi)\in
Hom_{\tilde{\cf}_{\lambda_{r}}^{\lambda_{l}}}(\Phi(A),\Phi(B))$
and, of course, the value of the functor morphism
$\Phi(\phi)$ on $V_{\lambda_{r}}^{k}$ corresponds to
$T^{\lambda_{l}}_{\mu}(\phi)\in Hom_{\hgtg}(\Phi(A)(V_{\lambda_{r}}^{k}),
\Phi(B)(V_{\lambda_{r}}^{k}))$:

\[i(\Phi(\phi))=T^{\lambda_{l}}_{\mu}(\phi).\; \qed\]
\end{sloppypar}

\bigskip

{\bf Definition.}
Let $(\lambda_{r},\lambda,\lambda_{l})$ be a translation
datum.
 Define the {\em restricted}
affine Harish-Chandra category $\widetilde{Hch}(\lambda_{l}
,\lambda_{r})$ to be the complete subcategory
$\Phi(\tco^{\lambda})\subset
 \tilde{\cf}_{\lambda_{r}}^{\lambda_{l}}$.

\bigskip

Conjecture \ref{conjhowdefhch} implies that
$\widetilde{Hch}(\lambda_{l}
,\lambda_{r})$ is equivalent to $\tco^{\lambda_{l}}$
and, in particular, independent of $\lambda$. This all
is in precise analogy with the Bernstein-Gelfand
theorem, see equivalence (\ref{bgisocat1}) in
Introduction. As a corollary, we get that the simple
objects of $\widetilde{Hch}(\lambda_{l}
,\lambda_{r})$ are in one-to-one correspondence with
the subset of the affine Weyl group $W_{k}$:
\[\{w\in W_{k}: w\cdot \lambda_{l}\in P^{+}\}.\]
Similarly, the functors $p_{\lambda_{l}}\circ(V_{w\cdot
\lambda}^{k}\tens ?)$ are obvious analogues of
the principal series representations.

A drawback of our definition is that
$\widetilde{Hch}(\lambda_{l}
,\lambda_{r})$ is defined only if there is $\lambda$ such
that $(\lambda_{r},\lambda,\lambda_{l})$ is a translation
datum. The simplest example when
$\widetilde{Hch}(\lambda_{l},\lambda_{r})$ is not defined is when
$\gtg=\gtsl_{2}, k=2,\lambda_{l}=1,\lambda_{r}=2$.
This drawback, however, is not as serious as it may seem
to be. In the case of special interest $\lambda_{l}=
\lambda_{r}$, the triple $(\lambda_{l},0,\lambda_{r})$
is a translation datum.

\subsection{Non-restricted case.}
One would prefer to have as many simple  objects as there are
elements in the entire affine Weyl group. To achieve that we use the results of
sect.\ref{genfromtcotook}.

Just as it was above, for a translation datum
$(\lambda_{l},\lambda,\lambda_{r})$
we have the category of functors $\cf_{\lambda_{r}}^{\lambda_{l}}$ from
$\tco^{\lambda_{r}}$ to $\co_{k}^{\lambda_{l}}$ and the  functor
\[\Phi:\co_{k}^{\lambda}\rightarrow \cf_{\lambda_{r}}^{\lambda_{l}},\]
\[ \Phi(A): B\mapsto p_{\lambda_{l}}(A\tens B).\]

There again arises the natural map

\[i:\; Hom_{\cf_{\lambda_{r}}^{\lambda_{l}}}(F,G)
\rightarrow Hom_{\hgtg}(F(V_{\lambda_{r}}^{k}),G(V_{\lambda_{r}}^{k})),\]
where $i(\psi)$ is  the value of the functor morphism $\psi$
on $V_{\lambda_{r}}^{k}$.

\begin{conjecture}
\label{conjhowdefhchnonr}
If $(\lambda_{r},\mu,\lambda_{l})$ and
 $(\lambda_{r},\nu,\lambda_{l})$ are translation data and
$A\in\co^{\mu}_{k},\;B\in\co^{\nu}_{k}$, then the map
\[Hom_{\cf_{\lambda_{r}}^{\lambda_{l}}}(\Phi(A),\Phi(B))
\rightarrow Hom_{\hgtg}(\Phi(A)(V_{\lambda_{r}}^{k}),
\Phi(B)(V_{\lambda_{r}}^{k})).
\]
is an isomorphism.
\end{conjecture}

Surjectivity of the map in Conjecture \ref{conjhowdefhchnonr}
in the case $\mu=\nu$ is proved just like surjectivity of the
map in Conjecture \ref{conjhowdefhch} except that instead of
Theorem
\ref{equivofcat}, one uses Theorem \ref{eqofcat_forbigg}.

We then define the Harish-Chandra category $HCh(\lambda_{l},\lambda_{r})$
as a complete subcategory of $\cf_{\lambda_{r}}^{\lambda_{l}}$ generated
by $\Phi(\co_{k}^{\lambda})$ if $(\lambda_{l},\lambda,\lambda_{r})$
is a translation datum. Provided Conjecture \ref{conjhowdefhchnonr}
is valid, this category is isomorphic to $\co_{k}^{\lambda_{l}}$.
Analogues of the principal series representations are, therefore,
$\Phi(M_{w\cdot\lambda}^{k}), w\in W_{k}$.

\end{document}